\documentclass[a4paper,  11pt, leqno]{amsart}

\usepackage[mathscr]{eucal}
\usepackage{amsmath}
\usepackage{amssymb}
\usepackage{amsfonts}
\usepackage{amsthm}

\usepackage{cases}

\usepackage[mathscr]{eucal}
\usepackage{amsmath}
\usepackage{amssymb}
\usepackage{amsfonts}
\usepackage{amsthm}
\usepackage[dvips]{graphics, color}
\theoremstyle{plain}
\newtheorem{theorem}{Theorem}[section]
\newtheorem{proposition}{Proposition}[section]
\newtheorem{remark}{Remark}[section]

\newtheorem{example}{Example}[section]
\newtheorem{definition}{Definition}[section]

\numberwithin{equation}{section}

\def\<{\left<} \def\>{\right>}

\def\proof{\noindent{\it Proof. }}
\def\bea{\begin{eqnarray} }
\def\eea{\end{eqnarray} }
\def\be{\begin{equation} }
\def\ee{\end{equation} }
\def\qed{\ifhmode\unskip\nobreak\fi\ifmmode\ifinner\else\hskip5pt
\fi\fi\hbox{\hskip5 pt \vrule width4 pt height6 pt depth1.5 pt \hskip1pt }}

\begin{document}
\title[]{A short note on biharmonic submanifolds in non-Sasakian 
contact metric $3$-manifolds}
\author[]{Toru Sasahara}
\address{Center for Liberal Arts and Sciences, 
Hachinohe Institute of Technology, 
Hachinohe, Aomori 031-8501, Japan}
\email{sasahara@hi-tech.ac.jp}


\begin{abstract}
{\footnotesize We characterize   biharmonic
 anti-invariant surfaces 
in  $3$-dimensional 
generalized $(\kappa, \mu)$-manifolds with  non-zero constant mean curvature
by means of the scalar curvature of the ambient space and the mean curvature.
In addition, we give a  method  for constructing infinity many examples of  
biharmonic submanifolds in a certain $3$-dimensional 
generalized $(\kappa, \mu)$-manifold.
Moreover, we  determine  $3$-dimensional  generalized $(\kappa, \mu)$-manifolds 
which admit a certain kind of proper biharmonic  foliation. 
}
\end{abstract}



\maketitle

\section{Introduction}
The notion of biharmonic submanifolds is a  natural extension of the notion of minimal submanifolds
 from a  variational point of view.
 Non-minimal biharmonic submanifolds are said to be proper.
Considerable advancement
 has been made in the study
of proper
biharmonic submanifolds 
 in manifolds with special metric properties (e.g., real space forms, complex space forms, Sasakian
space forms, conformally flat spaces, etc.) since the beginning of this century.

When the dimension of the ambient space is three, it seems worthwhile and 
interesting to 
construct and classify proper
biharmonic submanifolds in contact metric $3$-manifolds.
Some classification results for 
proper biharmonic submanifolds in Sasakian $3$-manifolds have been obtained
(see, for example, \cite{cad}, \cite{fc}, \cite{ino}).

As an extension of the notion of a Sasakian manifold, Koufogiorgos and 
Tsichlias \cite{kt} introduced the notion of   a generalized $(\kappa, \mu)$-manifold
 for two real
functions $\kappa$ $(\leq 1)$ and $\mu$.   If $\kappa\equiv 1$, then it is a Sasakian manifold.
Markellos and Papantoniou \cite{mp} studied  proper biharmonic Legendre curves and 
anti-invariant
surfaces of $3$-dimensional non-Sasakian generalized $(\kappa, \mu)$-manifolds.
However, they provided no examples of such submanifolds.
As far as the author know,
there have been no examples of proper biharmonic submanifolds
in non-Sasakian contact metric $3$-manifolds.

In this  short note,  
we first give a  method  for constructing infinity many examples of  proper
biharmonic Legendre curves in a certain non-Sasakian generalized $(\kappa, \mu)$-manifold.
This result corrects an error in  \cite{mp}. 
 We then give a characterization of    biharmonic
 anti-invariant surfaces 
in  $3$-dimensional non-Sasakian
generalized $(\kappa, \mu)$-manifolds with  non-zero constant mean curvature,
by means of the scalar curvature of the ambient space and the mean curvature.
Applying the result,  we obtain a  method  for constructing infinity many examples of   
proper biharmonic anti-invariant surfaces in a
 certain $3$-dimensional  non-Sasakian
generalized $(\kappa, \mu)$-manifolds.
Moreover, we determine 
$3$-dimensional   non-Sasakian generalized $(\kappa, \mu)$-manifolds which 
admit a certain kind of proper biharmonic anti-invariant foliation.

\section{Preliminaries}
Let $M^n$ be  
 an $n$-dimensional  submanifold of a Riemannian manifold $\tilde M$.
 Let us denote by $\nabla$ and
 $\tilde\nabla$ the Levi-Civita connections on $M^n$ and 
 $\tilde M$, respectively. The
 Gauss and Weingarten formulas are respectively given by
\bea 
 \begin{split}
 \tilde \nabla_XY&= \nabla_XY+B(X,Y), \label{gawe}\\
 \tilde\nabla_X N&= -A_{N}X+D_XN
 \end{split}
\eea 
 for tangent vector fields $X$, $Y$ and normal vector field $N$,
 where $B$, $A$ and $D$ are the second fundamental 
 form, the shape operator and the normal
 connection.

 The mean curvature vector field $H$ is defined by 
 \be H=(1/n){\rm trace}\hskip2pt B.\label{mean}\ee
The function $|H|$ is called the  {\it mean curvature}.
If it vanishes identically, then $M$ is called a {\it minimal submanifold}.

\section{Biharmonic submanifolds}
Let $f:M \rightarrow N$ be a smooth map of an
 $n$-dimensional Riemannian manifold
into another Riemannian manifold.
The {\it tension field}
$\tau(f)$ of $f$ is a section
of the induced vector bundle $f^{*}TN$
defined by
\be\tau(f)=
\sum_{i=1}^{n}\{\nabla^{f}_{e_i}df(e_i)
-df(\nabla_{e_i}e_i)\}\nonumber
\ee
for a local orthonormal frame $\{e_i\}$ on $M$, where  $\nabla^f$ and $\nabla$ denote
 the induced connection and the Levi-Civita connection of $M$, respectively.
If  $f$ is an isometric immersion,  then we have
\be\tau(f)=nH.\label{tau}
\ee

A smooth map $f$ is called
a {\it harmonic map} if 
it is a critical 
point of the energy functional
$$
E(f)=\int_{\Omega}|df|^2dv_g\nonumber 
$$
over every compact domain $\Omega$ of $M$, where $dv_g$ is the volume form of $M$. 
A smooth map $f$ is harmonic if and only if $\tau(f)=0$ at each point on $M$.

Eells and Lemaire introduced the notion of biharmonic maps as a natural generalization of 
the notion of  harmonic maps from a variational point of view. 
\begin{definition}[\cite{el}]
{\rm A smooth map $f: M\rightarrow \tilde M$ is  called a
{\it biharmonic map} if it is a critical point of the bienergy
functional defined by
\be
E_{2}(f)=\int_{\Omega}|\tau(f)|^2dv_g \nonumber
\ee
with respect to all variations with compact support.}
\end{definition}

The Euler-Lagrange equation for $E_{2}$ is given by 
\be \tau_{2}(f):=-\Delta_f(\tau(f))+{\rm trace}\tilde R(\tau(f),df)df=0, \label{eq:EL}
\ee
where
$\Delta_f=-{\rm trace}(\nabla^{f}\nabla^{f}-\nabla^{f}_{\nabla})$ and 
$\tilde R$ is 
the curvature tensor  of $\tilde M$  which is defined by $$\tilde R(X, Y)Z=[\tilde\nabla_X, \tilde\nabla_Y]Z-\tilde\nabla_{[X, Y]}Z$$ for the Levi-Civita connection $\tilde\nabla$ of $\tilde M$.

 If $f$ is a biharmonic isometric immersion, then $M$ or $f(M)$
 is called a {\it biharmonic submanifold} in $\tilde M$.  
It follows from (\ref{tau}) and (\ref{eq:EL})  that any
minimal submanifold is biharmonic. 
However,  the converse is not true in general.  Non-minimal biharmonic submanifolds are called
 {\it proper} biharmonic submanifolds.


 \section{Generalized $(\kappa, \mu)$-manifolds}
 A differentiable  manifold $\tilde M$ is called an {\it almost contact manifold}  if it
admits a unit vector field $\xi$, a one-form $\eta$ and a $(1, 1)$-tensor field $\phi$ satisfying
\be
\eta(\xi)=1, \quad 
\phi^2=-I+\eta\otimes\xi.\nonumber
\ee
Every almost contact manifold admits a Riemannian metric $g$ satisfying
$$g(\phi X, \phi Y)=g(X, Y)-\eta(X)\eta(Y).$$
If,  in addition,  the condition
$$d\eta(X, Y):=(1/2)(X(\eta(Y))-Y(\eta(X))-\eta([X, Y]))=g(X, \phi Y)$$ holds, then
$(\tilde M, \xi, \eta, \phi, g)$ is called a {\it contact metric manifold}.
A contact metric manifold  is called a 
{\it Sasakian manifold} if it satisfies
$$[\phi, \phi]+2d\eta\otimes \xi=0,$$ where
$[\phi, \phi]$ is the Nijenhuis torsion of $\phi$
which is defined as
$[\phi, \phi](X, Y)=\phi^2[X, Y]+[\phi(X), \phi(Y)]-\phi[\phi(X), Y]-\phi[X, \phi(Y)]$.

Koufogiorgos and Tsichlias \cite{kt} introduced the notion of 
{\it generalized $(\kappa, \mu)$-manifold},
which is defined as contact metric manifold whose curvature tensor $\tilde R$ satisfies 
\begin{equation}
\tilde R(X,Y)\xi =\left( \kappa I+\mu h\right) \left( \eta (Y)X-\eta (X)Y\right) \nonumber
\end{equation}%
for any vector field $X$ and $Y$, where
$2h$ is the Lie differentiation of $\phi$ with respect to $\xi$, 
and $\kappa$, $\mu$ are smooth functions. 
If $\kappa$ and $\mu$ are constant, then the manifold is simply called a {\it $(\kappa, \mu)$-manifold} (\cite{bkp}). 
Sasakian manifolds are $(\kappa, \mu)$-manifolds with $\kappa=1$ and $h=0$.

In \cite{kt}, it was proved that if the dimension of generalized $(\kappa, \mu)$-manifold is greater than $3$, then 
$\kappa$ and $\mu$ must be  constant, and if the dimension is $3$, then there exist examples of generalized $(\kappa, \mu)$-manifold
with $\kappa$, $\mu$ non-constant smooth functions.

Perrone \cite{pe} proved that a contact metric $3$-manifold $\tilde M^3$ is a generalized $(\kappa, \mu)$-manifold on an everywhere
dense open subset if and only if the vector field $\xi$
defines a harmonic map from $\tilde M^3$ into its unit tangent bundle equipped with
the Sasaki metric.

Any $3$-dimensional generalized $(\kappa, \mu)$-manifold $\tilde M^3$ with $\kappa<1$
 admits three mutually orthogonal distributions $D(0)={\rm Span}\{\xi\}$, 
$D(\lambda)$ and $D(-\lambda)$ determined by the eigenspaces of $h$, where $\lambda=\sqrt{1-\kappa}$. Furthermore, on such a manifold
  the following relations hold (see   Lemma 3.3 of \cite{kt1} and its proof): 
\bea
&& \tilde R(\phi X, X)X=\biggl(-\frac{\Delta\lambda}{2\lambda}-\frac{||{\rm grad}\hskip2pt\lambda||^2}{2\lambda^2}-\kappa-\mu\biggr)\phi X,\label{eq:ku4}\\
&& S=-\frac{\Delta\lambda}{\lambda}-\frac{||{\rm grad}\hskip2pt\lambda||^2}{\lambda^2}+2(\kappa-\mu), \label{eq:ku8}
\eea 
where $X$ is a unit vector  lying in $D(\pm\lambda)$, 
$\Delta\lambda=-\sum_{i=1}^{3}\{e_i(e_i\lambda)-(\tilde\nabla_{e_i}{e_i})\lambda\}$,
and $S$ denotes the scalar curvature.
\



Below, we shall exhibit examples of $3$-dimensional 
generalized $(\kappa, \mu)$-manifolds which are not $(\kappa, \mu)$-manifolds.
\begin{example}[\cite{kt1}]\label{ex1}
{\rm Let $\lambda:I\subset{\mathbb R}\rightarrow{\mathbb R}$ be a non-constant {\it positive}
 function defined on an open interval $I$.
 We put 
  $\lambda^{\prime}=\frac{d\lambda}{dz}$.
 Consider a manifold $\tilde M^3=\{(x, y, z)\in {\mathbb R}^2\times I \subset{\mathbb R}^3\}$.
 The vector fields
 \be
e_1=\frac{\partial}{\partial x}, \quad e_2=\frac{\partial}{\partial y},\quad
e_3=(\pm 2y+f(z))\frac{\partial}{\partial x}+\Bigl(2\lambda x-\frac{\lambda^{\prime}}{2\lambda}y
  +h(z)\Bigr)\frac{\partial}{\partial y}+\frac{\partial}{\partial z}\label{vec}
 \ee
are linearly independent at each point of $\tilde M^3$, where
$f(z)$ and $h(z)$ are arbitrary functions of $z$.
Let $g$ be the Riemannian metric defined by $g(e_i, e_j)=\delta_{ij}$, $i, j=1, 2, 3$,
and 
$\eta$ the dual $1$ form to $e_1$. 
If we define the $(1, 1)$-tensor field $\phi $ by $\phi e_1=0$, $\phi e_2=\pm e_3$ and $\phi e_3=\mp e_2$, then 
the manifold $(\tilde M^3, \phi , e_1, \eta ,g)$ 
  is a generalized $(\kappa, \mu)$-manifold
with
\be
\kappa=1-\lambda^2, \quad \mu=2(1\pm\lambda), 
\label{eq:ku}
\ee
where  double signs correspond to $\pm 2y$ in (\ref{vec}).  Let us 
 denote the manifold by $\tilde M(1-\lambda^2, 2(1\pm\lambda))$.
}
\end{example}

The following classification result was obtained by Koufogiorgos and Tsichlias.
 \begin{theorem}[\cite{kt1}]\label{cla1}
Let $\tilde M^3$  be a $3$-dimensional generalized $(\kappa, \mu)$-manifold which is not  a
$(\kappa, \mu)$-manifold. If $\tilde M^3$ satisfies  $\xi\mu=0$, then it 
 is locally given by
$\tilde M(1-\lambda^2, 2(1\pm\lambda))$.
\end{theorem}

\begin{remark}\label{rm1}
{\rm Every generalized  $(\kappa, \mu)$-manifold satisfies $\xi\kappa=0$ (see \cite{kt2}). }
\end{remark}

\begin{remark}{\rm
 Hereafter, except (\ref{anti-h}), all double signs correspond to that of $\tilde M(1-\lambda^2, 2(1\pm\lambda))$.}
 \end{remark}


On $\tilde M(1-\lambda^2, 2(1\pm\lambda))$, 
the following relations are true (see the proof of Theorem 4.2 in \cite{kt1}):
\bea
&&\tilde\nabla_{e_1}e_1=\tilde\nabla_{e_3}e_3=0,\label{eq:ku6}\\
&&\tilde\nabla_{e_2}e_2=\frac{\lambda^{\prime}}{2\lambda}e_3,\label{eq:ku7}\\
&&\tilde\nabla_{e_2}e_3=-\frac{\lambda^{\prime}}{2\lambda}e_2+(\lambda\pm 1)e_1,\label{eq:ku2}\\
&& he_2=\pm\lambda e_2. \label{eq:ku9}
\eea


\section{Biharmonic Legendre curves 
in $3$-dimensional generalized $(\kappa, \mu)$-manifolds}

Let $\tilde M^3$ be a contact metric $3$-manifold.
A  curve $\gamma: I\subset \mathbb{R}\rightarrow \tilde M^3$ parametrized by arclength
is called a {\it Legendre curve} if  $\eta(\gamma^{\prime}(s))=0$ holds identically.

For Legendre curves in $3$-dimensional $(\kappa, \mu)$-manifolds, we have
\begin{theorem}[\cite{mp}]
Let $\gamma$ be a non-geodesic Legendre curve in  a $3$-dimensional $(\kappa, \mu)$-manifold with
$\kappa<1$. Then $\gamma$ is biharmonic if and only if
either $\nabla_{\gamma^{\prime}}\gamma^{\prime}\parallel \phi\gamma^{\prime}$
and in this case $\gamma$ is a helix satisfying $k_g^2+\tau_g^2=-(\kappa+\mu)$ or 
$\nabla_{\gamma^{\prime}}\gamma^{\prime}\parallel \xi$ and 
$\gamma$ is a helix satisfying  $k_g^2+\tau_g^2=\kappa+\delta\mu$, where
$k_g$ {\rm (}resp. $\tau_g${\rm )} denotes the geodesic curvature {\rm(}resp. the geodesic torsion{\rm )}, and 
$\delta:=g(h\gamma^{\prime}, \gamma^{\prime})$ is constant.

\end{theorem}

On the other hand, in \cite[Theorem 3.2]{mp}  the authors claimed the following:

Let $\gamma$ be a Legendre curve satisfying 
$\nabla_{\gamma^{\prime}}\gamma^{\prime}\parallel \phi\gamma^{\prime}$
in a $3$-dimensional generalized $(\kappa, \mu)$-manifold which is not  a $(\kappa, \mu)$-manifold.
Then, $\gamma$ is biharmonic if and only if it is a geodesic.

However, {\bf this claim is incorrect}. In fact, there exist infinity many non-geodesic biharmonic Legendre curves
in $\tilde M(1-\lambda^2, 2(1\pm\lambda))$, as shown below.

For  two constants $b$ and  $c$, we consider the Legendre
 curve in $\tilde M(1-\lambda^2, 2(1\pm\lambda))$ given by
\be
\gamma_{\{b, c\}}(s)=(b, s, c). \label{legen}
\ee
It follows from  (\ref{eq:ku7}) that the curve  (\ref{legen})
 satisfies 
$\nabla_{\gamma^{\prime}}\gamma^{\prime}\parallel \phi\gamma^{\prime}$.

\begin{proposition}\label{pro:sa}
 A Legendre curve $\gamma(s)$ 
 in   $\tilde M(1-\lambda^2, 2(1\pm\lambda))$ given by $\gamma_{\{b, c\}}(s)=(b, s, c)$
 is proper biharmonic 
 if and only if $c$ satisfies $\lambda^{\prime}(c)\ne 0$ and 
 \be
 \left.\lambda\lambda^{\prime\prime}-2(\lambda^{\prime})^2-8\lambda^2(1\pm\lambda)\right|_{z=c}=0.\label{legen2}
\ee
\end{proposition}
\proof
The unit tangent vector field of  $\gamma$ is given by $e_2$. 
Equation (\ref{eq:ku7}) yields that 
 $\gamma$ is a geodesic if and only if $c$ satisfies
$\lambda^{\prime}(c)=0$. 
Equation  (\ref{eq:EL}) can be written as 
\be
\tilde\nabla_{e_2}\tilde\nabla_{e_2}\tilde\nabla_{e_2}e_2+\tilde R(\tilde\nabla_{e_2}e_2, e_2)e_2=0. \label{EL2}
\ee

Since (\ref{eq:ku9}) is satisfied, we put $X=e_2$ in (\ref{eq:ku4}).
Then, using (\ref{eq:ku6})-(\ref{eq:ku2}) and (\ref{eq:ku4}), by  a straightforward computation we can show 
that (\ref{EL2}) is equivalent to
(\ref{legen2}).  This finishes the proof.
\qed
\vspace{1.0ex}

\begin{example}\label{ex}{\rm 
If we put $\lambda=z^{-n}$ for $z>0$ and $n\in (0, 1)$, then $\lambda>0$, $\lambda^{\prime}\ne 0$ and 
(\ref{legen2}) becomes
\bea
8c^2(1\pm 
c^{-n})=n(1-n).\nonumber
\eea
For any $n\in (0, 1)$, this equation for $c$  has a positive solution, because $\lim_{c\rightarrow 0}c^2(1\pm 
c^{-n})=0$ and $\lim_{c\rightarrow\infty}c^2(1\pm 
c^{-n})=\infty$. }
 \end{example}

\section{Biharmonic anti-invariant surfaces in   $3$-dimensional generalized $(\kappa, \mu)$-manifolds}

A submanifold $M$ in a contact metric manifold is said to be {\it anti-invariant} if it is tangent to $\xi$ and satisfies 
$g(X, \phi Y)=0$ for all tangent vectors $X$ and  $Y$ of $M$.
For proper biharmonic anti-invariant surfaces in $3$-dimensional $(\kappa, \mu)$-manifolds, we have the following result.
\begin{theorem}[\cite{sasa}]
Let $f:M^2\rightarrow \tilde M^3$ be an anti-invariant isometric 
immersion of a  Riemannian $2$-manifold into a 
$3$-dimensional $(\kappa, \mu)$-manifold.
Then $f$ is proper biharmonic if and only if 
$\kappa=1 ;$ that is, $\tilde M^3$ is a Sasakian manifold, and moreover, 
$|H|^2=(1/8)(S-6)={\rm constant}$ $(\ne 0)$ on $M^2$.
\end{theorem}

This motivates us to study
 proper biharmonic anti-invariant surfaces 
in  $3$-dimensional 
generalized $(\kappa, \mu)$-manifolds which are not $(\kappa, \mu)$-manifolds.
The following theorem characterizes such surfaces with constant mean curvature by means of 
$S$ and $H$.
\begin{theorem}\label{chara}
Let $f: M^2\rightarrow \tilde M^3$ be an anti-invariant immersion of a Riemannian
$2$-manifold into a 
$3$-dimensional 
generalized $(\kappa, \mu)$-manifold which is not a  $(\kappa, \mu)$-manifold. 
Assume that $M^2$ has   non-zero constant mean curvature.
Then, $f$ is  biharmonic if and only if $S$ is constant on $M^2$
and
\bea
h(\phi H)=\biggl(\pm\sqrt{\frac{S}{6}-\frac{4}{3}|H|^2}-1\biggr)\phi H
\label{anti-h}
\eea
holds at any point of $M^2$.
\end{theorem}
\proof
Let  $M^2$ be an anti-invariant surface in a $3$-dimensional 
generalized $(\kappa, \mu)$-manifold which is not a  $(\kappa, \mu)$-manifold.
Assume that $M^2$ has  non-zero constant
mean curvature $\alpha$.  Let 
  $\{\xi, e_1, \phi e_1\}$ be a local orthonormal frame of $M^2$
 such that $H=\alpha\phi e_1$.
 We put \be
 \beta=1+g(he_1, e_1), \quad \gamma=g(he_1, \phi e_1)\label{anti-beta}.\ee
 By the same way as in the proof of Theorem 4.3 in \cite{mp} and using (4.9) of Proposition 4.1, 
 we see that
$M^2$ is  biharmonic if and only if $\beta$
 is constant on $M^2$,  and moreover, the following 
relations hold:
 \be
 \gamma=0,\quad
 -4\alpha^2-2\beta^2+(S/2)-\kappa-\mu(\beta-1)=0.\label{anti2}
 \ee

For any anti-invariant surface in  a 
$3$-dimensional 
generalized $(\kappa, \mu)$-manifold, 
the following relations are true
 (see Lemma 4.1 of \cite{mp}):
\bea
 &&(\beta-1)^2+\gamma^2=1-\kappa.\label{anti8}\\
 &&\xi(\gamma)=(2\beta-\mu)(\beta-1).\label{anti9}
\eea
Since $\kappa<1$,  it follows from (\ref{anti8}) and (\ref{anti9}) that
 the first equation in (\ref{anti2}) implies 
\be
\kappa=2\beta-\beta^2, \quad \mu=2\beta. \label{anti5}
\ee
By (\ref{anti-beta}), we see  that the first equation in (\ref{anti2}) is  equivalent to 
\be
h(\phi H)=(\beta-1)\phi H. \label{anti7}
\ee
Substituting (\ref{anti5}) into the second equation of (\ref{anti2}), we obtain
\be
-4\alpha^2-3\beta^2+(S/2)=0.\label{anti6}
\ee
Combining (\ref{anti7}) and (\ref{anti6}) leads to (\ref{anti-h}).
 Accordingly, we conclude that (\ref{anti2}) is equivalent to (\ref{anti-h}).
Under the condition of  (\ref{anti-h}), $\beta$ is constant if and only if $S$ is constant on $M^2$.
The proof is finished.
\qed
\vspace{1.0ex}

 Applying Theorem \ref{chara}, we obtain a method for constructing of infinity many 
examples of 
  proper biharmonic anti-invariant surfaces in $3$-dimensional generalized
$(\kappa, \mu)$-manifolds which are not $(\kappa, \mu)$-manifolds, as shown below.

For a constant $c$, we consider the immersion 
\be
f_c(x, y)=(x, y, c)\label{eq:c=}\nonumber
\ee of ${\mathbb R}^2$
into $\tilde M(1-\lambda^2, 2(1\pm\lambda))$.
Suppose that ${\mathbb R}^2$ is equipped with the induced metric $f^{*}_cg$.
Then, $\{f^{*}_ce_1, f^{*}_ce_2\}$ forms an orthonormal frame of $f_c$. 
It is clear that $f_c$ is an anti-invariant immersion
whose unit normal vector field $N$ is $\phi(f^{*}_ce_2)=f^{*}_ce_3$, and
 $\{f_c({\mathbb R}^2)| c\in I\}$ defines a codimension one foliation on  
$\tilde M(1-\lambda^2, 2(1\pm\lambda))$.

\begin{proposition}\label{pro:sa}
 The anti-invariant isometric immersion $f_c:{\mathbb R}^2\rightarrow
 \tilde M(1-\lambda^2, 2(1\pm\lambda))$ given by $f_c(x, y)=(x, y, c)$ is proper biharmonic 
 if and only if $c$ satisfies $\lambda^{\prime}(c)\ne 0$ and 
 \be
 \left.\lambda\lambda^{\prime\prime}-2(\lambda^{\prime})^2-8\lambda^2(1\pm\lambda)^2\right|_{z=c}=0.\label{eq:foli}
\ee
\end{proposition}
\proof
It follows from (\ref{eq:ku})-(\ref{eq:ku7}) that 
(\ref{eq:ku8}) can be rewritten as
\bea
S=
\frac{\lambda^{\prime\prime}}{\lambda}-\frac{3(\lambda^{\prime})^2}{4\lambda^2}-2(1\pm\lambda)^2. 
\label{eq:rho}
\eea
By  (\ref{gawe}), (\ref{mean}) and  (\ref{eq:ku7}) we have
 \be
|H|^2=\left.(\lambda^{\prime})^2/(16\lambda^2)\right|_{z=c}.
 \label{eq:A}
\ee
Thus, the surface $f_c({\mathbb R}^2)$ is minimal if and only if $\lambda^{\prime}(c)=0$.
By using   (\ref{eq:rho}),  (\ref{eq:A}) and (\ref{eq:ku9}),  we can show that
(\ref{anti-h}) is equivalent to (\ref{eq:foli}).  The proof is finished.
\qed
\vspace{1.0ex}

\begin{example}{\rm 
If w put $\lambda=z^{-n}$ for $z>0$ and $n\in (0, 1)$, then $\lambda>0$, $\lambda^{\prime}\ne 0$ and 
 (\ref{eq:foli}) becomes
\be
8c^2(1\pm 
c^{-n})^2=n(1-n).\nonumber
\ee
Similarly to   Example \ref{ex}, for any $n\in (0, 1)$ this equation for $c$ has a positive solution.}
\end{example}

\begin{remark}
{\rm A generalized $(\kappa, \mu)$-manifold with $||{\rm grad}\hskip2pt \kappa ||=a$ (constant)$\ne 0$ is locally
given by $\tilde M(1-\lambda^2, 2(1+\lambda))$ or  $\tilde M(1-\lambda^2, 2(1-\lambda))$, 
where
 $\lambda=\sqrt{1-az-b}$ for some constant $b$ (see \cite[Proposition 4.4]{kt} and \cite[Theorem 5]{kt2}).
  In this case, 
 since $\lambda$ does not satisfy (\ref{eq:foli}) for any $c$, there exists no anti-invariant surface $f_c({\mathbb R}^2)$
 which is  proper biharmonic. 
This implies that anti-invariant surfaces described in   
\cite[Example 4.1]{mp} is {\it not} proper biharmonic.
Hence, Proposition \ref{pro:sa} provides  the first examples of proper biharmonic surfaces
in non-Sasakian contact metric $3$-manifolds.}
\end{remark}

The following theorem provides us  with many examples of $3$-dimensional 
generalized $(\kappa, \mu)$-manifolds foliated by proper 
biharmonic  surfaces.

\begin{theorem}
Let $\tilde M^3$ be a 
$3$-dimensional 
generalized $(\kappa, \mu)$-manifold which is not a $(\kappa, \mu)$-manifold.
If $\tilde M^3$ admits an integrable distribution  $D(0)\oplus D(\pm\lambda)$ whose leaves are proper
 biharmonic surfaces, then $\tilde M^3$ is locally given by  $\tilde M(1-\lambda^2, 2(1\pm\lambda))$, where 
$\lambda$ is a positive solution of 
\be
(\lambda^{\prime})^2=\beta\lambda^4+16\lambda^4 \ln\lambda\mp 32\lambda^3-8\lambda^2\ne 0
\label{lambda}
\ee
for some constant $\beta$.
\end{theorem}
\proof
Let $\{\xi , X, \phi X\}$ be a local orthonormal frame such that
\be
hX=\lambda X, \quad h\phi X=-\lambda\phi X.\nonumber
\ee
Then, 
by Lemma 3.1 and 3.2 of \cite{kt1} we have
\be
\quad [\xi, X]=(1+\lambda-\mu/2)\phi X, \quad [\xi, \phi X]=(\lambda-1+\mu/2)X.
 \label{bla}\nonumber
\ee
  Hence, 
  $D(0)\oplus D(1+\lambda)$ or $D(0)\oplus D(1-\lambda)$ is integrable if and only if 
$\mu=2(1+\lambda)$ or $\mu=2(1-\lambda)$, respectively.
By applying Theorem \ref{cla1} and Remark \ref{rm1},   
  we obtain that
 if $\tilde M^3$ admits an integrable distribution  $D(0)\oplus D(\lambda)$ or  $D(0)\oplus D(-\lambda)$, then
   $\tilde M^3$ is locally given by 
 $\tilde M(1-\lambda^2, 2(1+\lambda))$ or $\tilde M(1-\lambda^2, 2(1-\lambda))$, respectively.
 In this case,  each leaf  is given by $f_c$.
All leaves of a foliation $\{f_c({\mathbb R}^2)| c\in I\}$
are  biharmonic surfaces if and only if (\ref{eq:foli})  is satisfied for all
$c\in I$.  By solving the ODE,  we have (\ref{lambda}).
\qed

\end{document}